\documentclass[letterpaper, 10 pt, conference]{ieeeconf}
\usepackage{graphicx}
\usepackage{subfiles}
\usepackage{amsmath}
\usepackage{amssymb}
\usepackage{cases}
\usepackage{subcaption}
\usepackage{empheq}
\usepackage{url}
\usepackage[dvipsnames]{xcolor}

\newcommand\MayBeColored[2]{#2}
\newcommand\UpdateA[1]{\MayBeColored{BrickRed}{#1}}
\newcommand\UpdateB[1]{\MayBeColored{blue}{#1}}

\newcommand\Fmin{F_\mathrm{min}}
\newcommand\SlidingStartTime{t^\mathrm S_0}
\newcommand\MuProjection{\mathrm{Proj}_{\mu}}

\newcommand\RR{\mathbb{R}{}}
\newcommand\sgn{\operatorname{sgn}}

\newtheorem{remark}{Remark}
\newtheorem{assumption}{Assumption}
\newtheorem{theorem}{Theorem}
\newtheorem{lemma}{Lemma}

\IEEEoverridecommandlockouts
\title{\UpdateA{Sign Primal-Dual Gradient Dynamics for Supply-Demand Coordination\\under Binary Human Feedback}}

\author{Yuya Miyaoka$^{1}$,  Kengo Urata$^{2}$, and Masaki Inoue$^{1}$
\thanks{*This work was supported by Grant-in-Aid for Scientific Research (B), No.~25K01254 from JSPS.}%
\thanks{$^{1}$Yuya Miyaoka and Masaki Inoue are with the Department of Applied Physics and Physico-Informatics in the Faculty of Science and Technology, Keio University, 3-14-1 Hiyoshi, Kohoku-ku, Yokohama, Kanagawa, Japan
        {\tt\small miyaoka.yuya@keio.jp; minoue.z6@keio.jp}}%
\thanks{$^{2}$Kengo Urata is with NTT Network Service Systems Laboratories, NTT Inc., Tokyo 180-8585, Japan
        {\tt\small kengo.urata@ntt.com}}%
}

\begin{document}
\maketitle
\thispagestyle{empty}
\pagestyle{empty}

\begin{abstract}
    \UpdateA{This paper addresses supply-demand coordination, in which a supplier sets prices to balance the demands of multiple buyers against limited resources.
    Existing studies typically solve the supply-demand coordination problem with a distributed algorithm called primal-dual gradient dynamics (PDGD), assuming that each buyer discloses the value of its own utility or marginal utility function.
    In contrast, we assume that each buyer discloses only binary information, namely whether it wants more or less of the current amount, and propose Sign PDGD. Sign PDGD is a PDGD algorithm in which the buyers' update depends only on the sign information.}
    We analyze the convergence of Sign PDGD by interpreting it as a sliding mode control (SMC). For a quadratic cost and linear constraints, we show that under a \UpdateA{gain condition}, the \UpdateA{primal and dual variables} converge to the optimal solution.
    Numerical simulations on a toy resource-sharing problem confirm the two-phase behavior of SMC and show that the \UpdateA{primal and dual variables} approach the optimal solution.
\end{abstract}

\section{Introduction}

Sharing limited resources such as networks, electricity, and water among multiple buyers has become an important societal problem. In such settings, price-based mechanisms that balance buyer-side demand with supply-side capacity have long been used. The mechanism is called supply-demand coordination \cite{Siano14_SupplyDemandCoordination}\cite{Palomar06_NUM}.

Supply-demand coordination is typically modeled as a distributed gradient-based dynamics over the primal and dual variables of an optimization problem.
This form arises naturally from the roles of the buyer side and the supply side. Each buyer adjusts its purchase to maximize its own utility, while the supplier sets prices to enforce the finite-resource constraint. Each side carries out its role independently, and as a result, the distributed gradient-based dynamics emerge.

Such price-based, iterative coordination has been studied in various application domains.
In power systems, distributed schemes are developed for voltage determination in electric power networks \cite{Wang_2020_AsynchronousVoltage}, for dynamic price adjustment that eliminates power imbalance \cite{Hansen2014DynamicMarket}, and for real-time pricing that maximizes the users' utility through the interaction between users and the provider \cite{Samadi2010_OptimalRealTimePricing}.
In water resource management, allocation that maximizes the users' utility is presented \cite{Janjua_2025_ThreeStageWater}.
For computing resources, price adjustment that maximizes the users' utility has been formulated based on user-provider interaction and solved by the Alternating Direction Method of Multipliers (ADMM) \cite{Nguyen2019MarketBasedFog}. The work \cite{Moro_2023_JointComputeRadio} applies the market equilibrium theory for partitioning of computer resources.
From a theoretical viewpoint, the works \cite{Palomar06_NUM} and \cite{Feijer_2010_PrimalDualStability} interpret supply-demand coordination for network resources within the mathematical framework of primal and dual decomposition, and analyze the stability of the resulting primal-dual dynamics.
Their studies share a common structure: a loop in which the buyer side and the supplier side repeatedly exchange amounts and prices to reach a coordinated allocation.

However, most of these approaches implicitly assume that each buyer can quantitatively evaluate their own marginal utility or utility at any given resource amount. In other words, the quantity demanded is treated as a primal variable that the buyer adjusts continuously along the gradient of its utility. If the buyers are actual humans, this assumption is challenging to justify: a human buyer cannot know their own utility function precisely, nor can they express the magnitude of their marginal utility or change their demand smoothly and faithfully according to it. What a human buyer can realistically provide is much coarser information, such as whether they want to increase or decrease their current resource amount.

\begin{figure}
    \centering
    \includegraphics[width=0.9\linewidth]{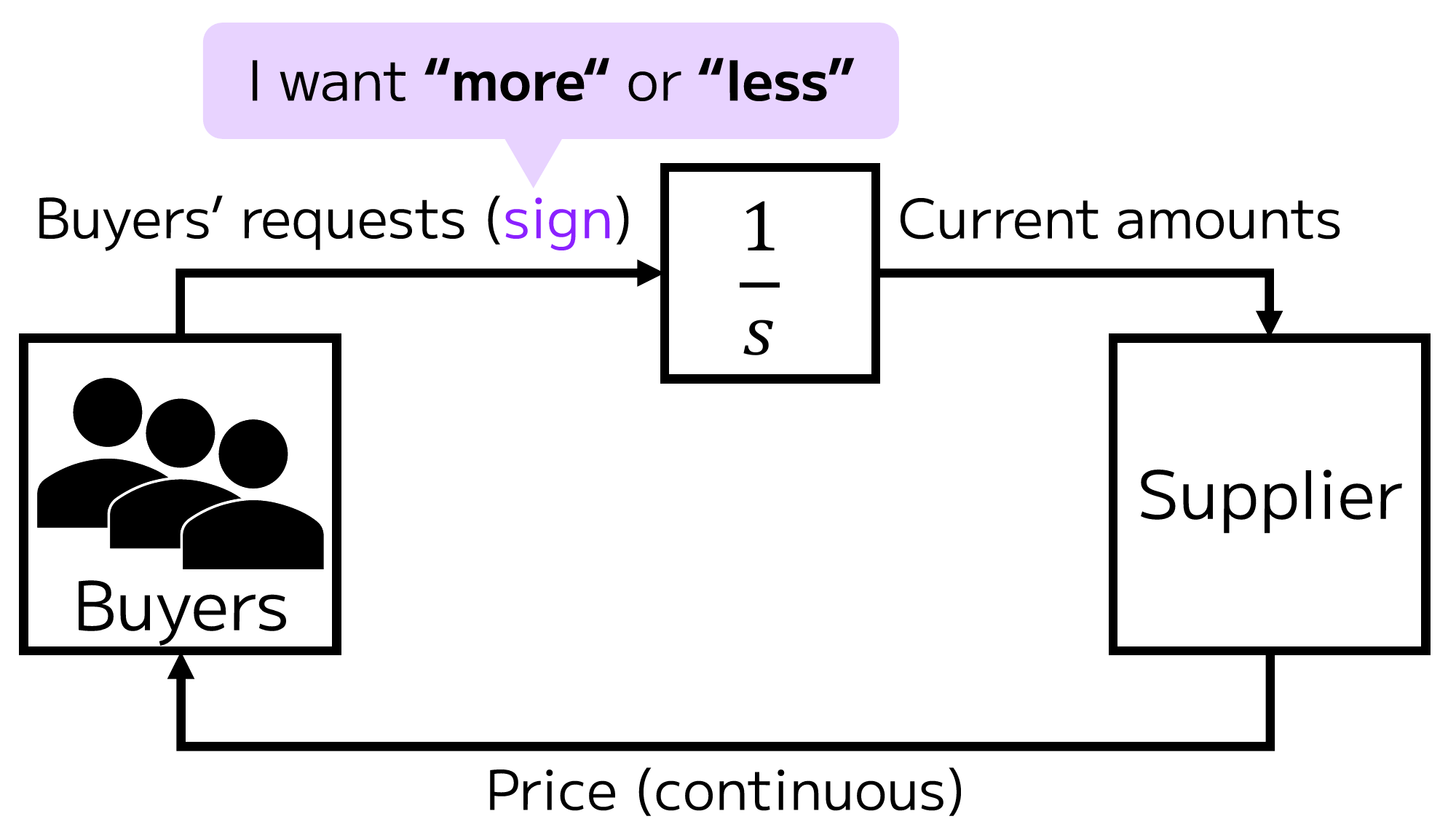}
    \caption{\UpdateA{The concept of the supply-demand coordination with human buyer feedback}}
    \label{F:I.SystemConcept}
\end{figure}

In this paper, \UpdateA{we introduce a more realistic human buyer model that provides only coarse information.}
Fig.~\ref{F:I.SystemConcept} illustrates the supply-demand coordination considered in this paper. 
During the coordination process, the resource allocation amount is adjusted according to the buyers' preferences, while the prices are adjusted by the supplier to satisfy resource constraints. The amount and the price are then updated repeatedly.
The distinctive point in this paper is that the buyer's means of expressing their intentions is narrow: at any given moment, a buyer is allowed to provide binary information, that is, whether they want ``more'' of the current amount or ``less'' of it. This motivates us to formulate a mathematical model in which information from the human buyer is available only as a sign signal.


\UpdateA{We propose sign primal-dual gradient dynamics (Sign PDGD)}, a PDGD in which the primal variable increases or decreases at a constant rate. We also analyze the convergence of Sign PDGD. \UpdateA{The key point of the analysis is that we incorporate the concept of sliding mode control (SMC): we divide the behavior of Sign PDGD into a reaching phase, in which the state reaches a certain manifold, and a sliding phase, in which the state moves along the manifold toward the optimum. We show that Sign PDGD converges to the optimum under a gain condition on the primal value update.}
Finally, we address a toy model of a supply-demand coordination problem as a numerical example and apply Sign PDGD.

Regarding the theoretical analysis, the convergence properties of general PDGD have been studied in \UpdateA{several works \cite{Yamashita20_PDGD}\cite{Dean15_PDGD}\cite{Sulaiman20_PDGD}\cite{Qu19_PDGD}\cite{Ashish16_PDGD}}.
However, to our knowledge, no existing work discusses the behavior and convergence of PDGD in which some of the variables evolve along a sign gradient, as in our Sign PDGD. The analysis in this paper fills this gap.

\section{Sign PDGD}
In this section, we formulate the supply-demand coordination problem and its \UpdateA{corresponding} primal-dual gradient dynamics (PDGD). Then, we modify the buyer's model to reflect coarse human feedback, yielding our proposed dynamics, Sign PDGD.

\subsection{Problem Setting} 
A general supply-demand coordination problem is formulated as follows:
\begin{subequations}
\begin{empheq}[left={\empheqlbrace}]{alignat=1}
    \min_x \quad& f(x), \label{E:M.Prob.f} \\
    \text{s.t.} \quad
        & g(x) = 0, \label{E:M.Prob.g} \\
        & c(x) \le 0, \label{E:M.Prob.c}
\end{empheq}\label{E:M.Prob}
\end{subequations}

\noindent \UpdateA{where} $x \in \RR^N$ denote the amounts purchased by the buyers, and $f: \RR^N \to \RR$ is a cost function representing the buyers' dissatisfaction\UpdateA{, i.e., the negative of the buyers' utility}. The smaller the value of $f$, the greater the satisfaction.  $g:\RR^N \to \RR^M$ and $c:\RR^N \to \RR^L$ are the supplier-side equality and inequality constraints, respectively.
When there are multiple buyers, $f$ is separable across the buyers. Notably, $-\nabla f(x)$ is the buyers' marginal utility: the rate at which satisfaction increases when a buyer acquires more of the resources. 
An example of $g$ is a constraint that requires the sum of the component ratios in the resource blending problem to equal $1$. An example of $c$ is a constraint that the \UpdateB{total resource demand does not exceed the available supply}.
\begin{remark} 
    Problem \eqref{E:M.Prob} \UpdateA{represents} a setting in which several suppliers provide several types of resources, and these resources are divided among several buyers. For example, $x$ can be arranged as \UpdateA{$x=(x_{11},\dots,x_{LI})^\top$, where $x_{li}$ is the amount of resource $l$} purchased by buyer $i$. The constraints can be imposed per supplier or per resource type, depending on the situation.
\end{remark}

\UpdateA{We let the Lagrangian of \eqref{E:M.Prob} be}
\begin{align}
    L(x, \lambda, \mu) = f(x) + \lambda^\top g(x) + \mu^\top c(x),
\end{align}\label{E:M.Lagrangian}

\noindent \UpdateA{where} $\lambda \in \RR^M$ and $\mu \in \RR^L$ are the Lagrange multipliers of \eqref{E:M.Prob.g} and \eqref{E:M.Prob.c}, respectively.
In the general supply-demand coordination problem, the dual variables play the role of prices that the supplier updates \cite{Palomar06_NUM}\cite{Samadi2010_OptimalRealTimePricing}\cite{Nguyen2019MarketBasedFog}.

\subsection{Conventional PDGD} 
The roles of the buyer and the supplier in supply-demand coordination are often expressed as a distributed process in the form of PDGD \cite{Palomar06_NUM}\cite{Feijer_2010_PrimalDualStability}\cite{Shomalzadeh_2025_SociallyAcceptable}. In PDGD, the buyer and supplier roles are clearly divided. The buyer updates the primal variable $x$ to maximize their own benefit, and the supplier updates the dual variables $(\lambda,\mu)$ to satisfy the constraints. 

The conventional continuous-time PDGD for \eqref{E:M.Prob} is given as follows:
\begin{subequations}
\begin{empheq}[left={\empheqlbrace}]{alignat=1}
    \dot x &= - \eta_x \nabla_x L(x, \lambda, \mu), \label{E:M.ConventionalPDGD.x} \\
    \dot \lambda &= \nabla_\lambda L(x, \lambda, \mu) = g(x), \label{E:M.ConventionalPDGD.lambda} \\
    \dot \mu &= [ \nabla_\mu L(x, \lambda, \mu) ]_\mu^+ = [c(x)]_\mu^+, \label{E:M.ConventionalPDGD.mu}
\end{empheq}\label{E:M.ConventionalPDGD}
\end{subequations}

\noindent \UpdateA{where} $\eta_x$ is \UpdateA{a positive constant representing a gain of the $x$ update}, and $[\cdot]_\mu^+$ is the projection that keeps $\mu$ nonnegative: for each element $j$,
\begin{align}
    ([c]_\mu^+)_j = \begin{cases}
        c_j, & \mu_j > 0 \ \text{or}\ c_j > 0, \\
        0, & \text{else}
    \end{cases}
\end{align}

\begin{remark} 
    \UpdateA{To clarify the meaning of the dynamics \eqref{E:M.ConventionalPDGD}}, consider linear constraints $g(x) = Ax + b$ and $c(x) = Cx + d$. Then, \eqref{E:M.ConventionalPDGD.x} is \UpdateA{reduced to}
    \begin{align}
        \nabla_x L = \nabla f(x) + A^\top \lambda + C^\top \mu,
    \end{align}
    where $(-\nabla f(x))$ is the marginal utility and $A^\top\lambda+C^\top\mu$ is the marginal payment determined by the current prices.
    In this sense, \eqref{E:M.ConventionalPDGD.x} represents that each buyer increases the acquired amount when the marginal utility exceeds the price, and decreases it otherwise.
    \eqref{E:M.ConventionalPDGD.lambda} and \eqref{E:M.ConventionalPDGD.mu} imply that the supplier raises the price of a constraint when it is violated, and lowers it otherwise.
\end{remark}

\UpdateA{A fundamental assumption underlying coordination via conventional PDGD \eqref{E:M.ConventionalPDGD} is that buyers have access to their own utility or marginal utility value. Although this assumption is reasonable when buyers are automated agents, it becomes impractical when buyers are human decision makers.} This motivates us to remodel the buyers' dynamics in the next subsection.

\subsection{Sign PDGD: Human Buyers with Coarse Feedback}
\UpdateA{We aim to update the dynamics of the primal variable $x$ in} the PDGD \eqref{E:M.ConventionalPDGD}.  We assume that at the current prices, a human buyer only expresses whether he/she wants \textit{more} or \textit{less} of the resource. That means that a buyer reports only the sign of $\nabla_x L$. The update logic of the primal variable \eqref{E:M.ConventionalPDGD.x} is replaced by $\dot x = -\eta_x \sgn(\nabla_x L(x,\lambda,\mu)) $, where $\sgn$ is an element-wise sign function. Note that the update logic of the dual variables, \eqref{E:M.ConventionalPDGD.lambda} and \eqref{E:M.ConventionalPDGD.mu}, are unchanged.

\UpdateA{The updated dynamics, called Sign PDGD, are described as}
\begin{subequations}
\begin{empheq}[left={\empheqlbrace}]{alignat=1}
     \dot x &= -\eta_x \sgn(\nabla_x L(x,\lambda,\mu)) , \label{E:M.SPDGD.x} \\
    \dot \lambda &= \nabla_\lambda L(x, \lambda, \mu) = g(x), \label{E:M.SPDGD.lambda} \\
    \dot \mu &= [ \nabla_\mu L(x, \lambda, \mu) ]_\mu^+ = [c(x)]_\mu^+. \label{E:M.SPDGD.mu}
\end{empheq}\label{E:M.SPDGD}
\end{subequations}

Sign-based gradient methods have also been studied in machine learning. A representative example is Sign Gradient Descent, which updates each parameter using the sign of its gradient. This approach avoids the vanishing gradient problem in training machine learning models \UpdateA{\cite{A23_SignGradient}\cite{Moulay19_SignGradient}\cite{Liao18_SignGradient}\cite{Riedmiller93_SignGradient}}.

\section{Theoretical Analysis}\label{Section:TheoreticalAnalysis}
\UpdateA{In this section, we address the convergence analysis of Sign PDGD} \eqref{E:M.SPDGD}. \UpdateA{A key feature of Sign PDGD is that the sign function is expected to induce finite-time convergence to a certain manifold. To exploit this property, we interpret the dynamics as a stabilization problem based on sliding mode control (SMC), which is} widely used for control applications, such as robotics \cite{Sen26_SMC}\cite{Mahmoud25_SMC}\cite{Liu20_SMC}, networks \cite{Mortajine25_SMC}, and power plants \cite{Surjagade22_SMC}.
\UpdateA{Then, we divide the convergence analysis into two phases. In the \textit{reaching phase}, the state approaches the manifold called the \textit{sliding surface}, and in the \textit{sliding phase},} the state remains on the sliding surface \cite{Filippov88_SMC}\cite{Weibing93_SMC}.

\subsection{Setting and Notation}
\UpdateA{Consider that $f$ is a separable quadratic objective funtion and $(c,g)$ are linear constraints.}
\begin{align}\begin{split}
    & f(x)=\frac12 (x-\bar x)^\top F(x-\bar x) ,\\
    & g(x)=Ax+b ,\quad c(x)=Cx+d,
\end{split}\end{align}
where $F=\mathrm{diag}(F_1,\dots,F_N)$ with $F_i>0$ for all $i\in\{1,2,\dots,N\}$, $\bar x\in\RR^N$, $A\in\RR^{M\times N}$, $C\in\RR^{L\times N}$, $b\in\RR^M$, and $d\in\RR^L$.
Let $\Fmin := \min_i F_i$.

The Lagrange multipliers $\lambda\in\RR^M$ and $\mu\in\RR^L$ are collected into $z:=(\lambda,\mu)$.
\UpdateA{We let}
\begin{align}\begin{split}
    G :=& [A^\top \ C^\top] ,\\
    H :=& G^\top F^{-1} G ,\\
    h:=& [(A\bar x+b)^\top \ (C\bar x+d)^\top]^\top .
\end{split}\end{align}
Then,
\begin{align}
    \nabla_x L =& F(x-\bar x)+v , \\
    v:=&\UpdateA{Gz}=A^\top\lambda+C^\top\mu.
\end{align}

We assume that $H\succ0$ and that there is no duality gap, so that an optimal primal solution $x^\ast$  and the dual solution $z^\ast = (\lambda^\ast,\mu^\ast)$ exist.
\UpdateA{From the stationarity condition of the Lagrangian with respect to $x$, we have} $F(x^\ast-\bar x)+Gz^\ast=0$.

\begin{assumption}[Price Velocity Limiter]\label{Assumption:A.PriceVelocityBound}
   For all $t\ge 0$ and for all $i$,
    \begin{align}
        \sup_{t\ge0}|\dot v_i(t)| < \eta_x F_i.
    \end{align}
\end{assumption}

\begin{remark}
    Assumption~\ref{Assumption:A.PriceVelocityBound} requires that the gain of the buyers' sign-based feedback, denoted by $\eta_x$, be fast enough to cancel the price drift.
    This is systematically achieved by applying a speed filter to \eqref{E:M.SPDGD.lambda} and \eqref{E:M.SPDGD.mu}.
\end{remark}

We let the variable $s$ be defined as follows:
\begin{align}\label{E:A.s}
    s := \nabla_x L = F(x-\bar x)+v .
\end{align}

\UpdateA{%
\noindent Sign PDGD first moves from $s\ne0$ toward $s=0$, and then switches to a motion within $s=0$, whose details are discussed in the following subsections. These two motions can be regarded as the phases of SMC. In this sense, $s$ is a \textit{sliding variable}, and the set $s=0$ is the \textit{sliding surface}.}

\subsection{Reaching Phase} 
In this subsection, we let $t=t^\mathrm R_0:=0$ be the initial time of the reaching phase, and also the initial time of the whole control time. Consider that $s(0)\ne0$, implying that the state variables $(x(0),\mu(0),\lambda(0))$ are not on the sliding surface yet.
We also let $s_0=s(0)$.

\begin{lemma}[Finite-Time Reaching]\label{L:A.Reach}
    Under Assumption~\ref{Assumption:A.PriceVelocityBound}, the sliding variable $s_i$ reaches zero no later than
    \begin{align}
       T^\mathrm R_i := \frac{|s_{0,i}|}{\eta_x F_i - \sup_{t\ge0}{|\dot v_i(t)|}},
       \label{E:A.R.ReachTime.0}
    \end{align}
    and all components have reached zero by
    \begin{align}
        T^\mathrm R := \max_i T^\mathrm R_i .
        \label{E:A.R.ReachTime.1}
    \end{align}
\end{lemma}

\begin{proof}
    Since $F$ is diagonal, differentiating \eqref{E:A.s} gives
    \begin{align}
        \dot s_i = F_i\dot x_i+\dot v_i = -\eta_xF_i \sgn(s_i)+\dot v_i
    \end{align}
    for each $i$. While in the reaching phase, i.e., $s_i\ne0$,
    \begin{align}
        \frac{d}{dt} |s_i|
        = \sgn(s_i)\dot s_i
        = -\eta_xF_i+\sgn(s_i)\dot v_i
        <0.
        \label{E:A.R.Decrease}
    \end{align}
    Hence $|s_i|$ decreases while it is nonzero, which gives the finite reaching time \eqref{E:A.R.ReachTime.0}.
\end{proof}

\subsection{Sliding Phase} 
In this subsection, we assume the reaching phase is complete for all components of $s$, and we set $t= \SlidingStartTime$ as the initial time of the sliding phase. Then, $s(\SlidingStartTime) = 0$ holds, implying $(x(\SlidingStartTime),\lambda(\SlidingStartTime),\mu(\SlidingStartTime))$ are on the sliding surface. 

The following lemma shows that the sliding phase sustains once it starts. 
\begin{lemma}[Sustaining the Sliding Phase]\label{L:A.S.Sustain2}
    Under Assumption~\ref{Assumption:A.PriceVelocityBound}, if $s(t_1)=0$ \UpdateA{at a time $t_1\ge0$}, then $s(t)=0$ for all $t\ge t_1$.
\end{lemma}
\begin{proof}
    We examine whether $s$ in \eqref{E:A.s} remains at $0$.
    For each $i$, differentiating \eqref{E:A.s} and using $F$ diagonal gives
    \begin{align}
        \dot s_i = F_i\dot x_i+\dot v_i=-\eta_x F_i\sgn(s_i)+\dot v_i .
    \end{align}
    On the sliding surface, the right-hand side is discontinuous, so the solution is understood in the sense of Filippov \cite{Filippov88_SMC}.
    Then, $\sgn(s_i)$ is replaced by a set-valued map taking values in $[-1,+1]$ at $s_i=0$.
    Hence, $s_i$ can stay at $0$ if there exists a value $\sigma_i\in[-1,+1]$ such that
   \begin{align}
       0=-\eta_xF_i\sigma_i+\dot v_i.
   \end{align}
   This equality holds if and only if
   \begin{align}
       \sigma_i=\frac{\dot v_i}{\eta_xF_i},
   \end{align}
   
   By Assumption~\ref{Assumption:A.PriceVelocityBound}, we have $-1<\sigma_i<1$, so it is an admissible Filippov selection.
   Therefore, $s_i(t)=0$ is sustained for all $t\ge t_1$, and since this holds for every $i$, the lemma holds.
\end{proof}

Now, we discuss the dynamics of the state $(x, z)$ in the sliding phase. 
By Lemma~\ref{L:A.S.Sustain2}, $s(t)=0$ for all $t\ge\SlidingStartTime$, and the primal variable $x$ follows the equivalent dynamics
\begin{align}
    x= \bar x-F^{-1} v = \bar x -F^{-1} Gz.
\end{align}
Substituting \UpdateA{$x$} into the dual variable dynamics in \eqref{E:M.SPDGD} gives
\begin{align}
    \dot \lambda &= A\UpdateA{x} + b = -AF^{-1}v + A\bar x+b, \\
    \dot \mu &= [C\UpdateA{x}+d]_\mu^+ = [-CF^{-1}v + C\bar x+d]_\mu^+.
\end{align}
Writing $u(z):=-Hz+h$, \UpdateA{we have}
\begin{align}\label{E:A.S.DualDynamics}
    \dot z = \MuProjection(u),
\end{align}
where $\MuProjection(z)$ is the projection function $[\cdot]_\mu^+$ only applying the $\mu$-part of $z$.
Let $I(t)$ denote the set of indices $i$ with $\dot z_i = u_i(z(t))$, meaning the projection $\MuProjection$ is inactive. This set always contains the indices of $\lambda$. Let $I^C(t)$ be its complement, which consists of the indices of $\mu$ where the projection is active.

The dual function of \eqref{E:M.Prob} is $D(z) = -\frac12 z^\top Hz + h^\top z$, and $u$ can be written using the dual gradient, $u = \nabla_z D$. \UpdateA{In this sense, the dynamics of $z$ in the sliding phase can be viewed as a dual ascent dynamics. The following lemma shows that $z$ with the dual ascent dynamics converges to the optimum.}

\begin{lemma}[Convergence of  the Dual Variable]\label{L:A.S.DualConv}
    Under the dual dynamics \eqref{E:A.S.DualDynamics}, for any initial value $z(\SlidingStartTime)=z_0$, the dual variable $z$ converges to the optimum $z^\ast$ as $t\to\infty$.
\end{lemma}
\begin{proof}
    Let $W(z):=\frac12\|z-z^\ast\|_2^2$. Along the dual dynamics \eqref{E:A.S.DualDynamics},  the function satisfies
    \begin{align}\begin{split}
        \dot W&= (z-z^\ast)^\top \dot z \\
        &= (z-z^\ast)^\top u - \sum_{i\in I^C} (z_i-z_i^\ast)u_i \\
        &= (z-z^\ast)^\top u + \sum_{i\in I^C} z_i^\ast u_i.
    \end{split}\end{align}
    \UpdateA{The second equality follows because $z_i=0$ for $i\in I^C$.}
    The last sum is nonpositive because $z_i^\ast\ge0$ and $u_i\le0$ for $i \in I^C$.

    Let $u=u(z)$ and $u^\ast=u(z^\ast)$. Since $u-u^\ast = -H(z-z^\ast)$,
    \begin{align}\begin{split}
        & (z-z^\ast)^\top u \\
        & = -(z-z^\ast)^\top H (z-z^\ast) + (z-z^\ast)^\top u^\ast.
    \end{split}\end{align}

    \UpdateA{Now, we divide $(z-z^\ast)^\top u^\ast$ into the $\lambda$-part and the $\mu$-part elements, as follows:
    \begin{align}\begin{split}
        &(z-z^\ast)^\top u^\ast \\
        &= \sum_{i\in(\lambda\text{-part})} (z_i-z_i^\ast)u_i^\ast
        + \sum_{i\in(\mu\text{-part})} (z_i-z_i^\ast)u_i^\ast
    \end{split}\end{align}
    }

    \noindent The KKT conditions for the dual maximization give \UpdateA{$u^\ast_i=0$ for $i$ in the $\lambda$-part, and $z^\ast_i \ u^\ast_i = 0$ for $i$ in the $\mu$-part.} Therefore
    \UpdateA{
    \begin{align}
        (z-z^\ast)^\top u^\ast
        = \sum_{i\in(\mu\text{-part})} z_i u_i^\ast .
    \end{align}
    Since $z_i\ge0$ and $u_i^\ast \le0$ for $i$ in the $\mu$-part, $(z-z^\ast)^\top u^\ast$ is nonpositive. 
    }
    
    Combining them yields 
    \begin{align}
        \dot W \le -(z-z^\ast)^\top H (z-z^\ast)
        \le - 2 \lambda_{\min}(H) W.
    \end{align}
\end{proof}

\UpdateA{From Lemmas~\ref{L:A.S.Sustain2} and \ref{L:A.S.DualConv}, we have the following theorem.}
\begin{theorem}[Convergence of Sliding Phase]\label{Theorem:A.S.Convergence}
    Consider that $s(\SlidingStartTime) = 0$, and we assume that the number of changes of $I(t)$ is finite.
    Then, the state variables $(x,\lambda,\mu)$ converge to the optimal values $(x^\ast,\lambda^\ast,\mu^\ast)$ as $t\to\infty$.
\end{theorem}
\begin{proof}
    By Lemma~\ref{L:A.S.Sustain2}, the sliding phase sustains from $t=\SlidingStartTime$ onward.
    During the sliding phase, the dual variable $z$ moves according to \eqref{E:A.S.DualDynamics}.
    Then, by Lemma~\ref{L:A.S.DualConv}, $z(t)\to z^\ast$, and
    \begin{align}\begin{split}
    & x(t) = \bar x - F^{-1}Gz(t) \\ 
    & \to \bar x-F^{-1}Gz^\ast = x^\ast,
    \end{split}\end{align}
    where the last equality follows from $Fx^\ast+Gz^\ast=0$.
\end{proof}

\begin{remark}
    The analysis is addressed in continuous time. In practice, supply-demand coordination proceeds in discrete time, and the discontinuous update of the Sign PDGD causes chattering around the sliding surface. Therefore, a numerical analysis for the discrete-time case, including the convergence and the size of the chattering, needs to be considered separately.
\end{remark}

\section{Numerical Simulation}\label{Section:Experiment3}
In this section, we apply the Sign PDGD \eqref{E:M.SPDGD} to a toy supply-demand coordination problem and compare their behavior with the theoretical analysis discussed in Section~\ref{Section:TheoreticalAnalysis}.

\subsection{Problem Setting}
We consider a single supplier and a single finite resource, such as the CPU capacity of a single datacenter, whose total stock is $V=7$. This resource is shared among $N=4$ buyers, and $x_i$ denotes the amount determined by buyer $i$.

The cost function of buyer $i$ is $f_i (x_i) = \frac12 F_i (x_i - \bar x_i)^2$,
where $\bar x_i>0$ is the amount that the buyer $i$ ideally wants, and $F_i>0$ represents the dissatisfaction factor of the buyer $i$. Note that the value of $\bar x_i$ is private, meaning that it is not open to other buyers or the supplier.
The parameters are set to $\bar x = [3.0,2.7,3.2,2.6]$ and $F=\mathrm{diag}(0.5,0.6,1.0,1.5)$. \UpdateB{We let the objective function as $f(x)=\frac12 (x-\bar x)^\top F (x-\bar x).$}

The only constraint is that the total amount must not exceed the stock $V$, which is formulated as
\begin{align}\label{E:E2.Constraint1}
    1_N^\top x- V\le 0.
\end{align}

\noindent There is no equality constraint in this setup.

\UpdateA{By letting $C=1_N^\top$ and $d=V$, we reduce the optimization problem \eqref{E:M.Prob} to}
\begin{subequations}
\begin{empheq}[left={\empheqlbrace}]{alignat=1}
    \min_{x\in\RR^4} \quad& f(x) \\
    \text{s.t.} \quad
        & \text{\eqref{E:E2.Constraint1}} .
\end{empheq}\label{E:E3.Prob}
\end{subequations}

\noindent The dual variable $\mu\in\RR$ is the scalar price of the resource, and $z=\mu$, $G=C^\top=1_N$, and $H=\sum_i 1/F_i$.

The analytical optimum is
\begin{align}
    \mu^\ast = 
    \frac{\sum_i\bar x_i-V}{\sum_i 1/F_i} \approx 0.84,
    \quad x_i^\ast = \bar x_i - \frac{\mu^\ast}{F_i},
\end{align}
which we use as the reference for the simulation results.

\subsection{Scenarios}
We compare two 
scenarios: the Sign PDGD \eqref{E:M.SPDGD} with $\eta_x=30.0$ and $3.0$. 
All the simulations are discretized with the time step $\Delta t = 0.001~\text{s}$, and the initial values are $x(0) = 0$ and $\mu(0) = 0$.

\subsection{Results for $\eta_x=30.0$}

\begin{figure}
    \centering
    \includegraphics[width=0.8\linewidth]{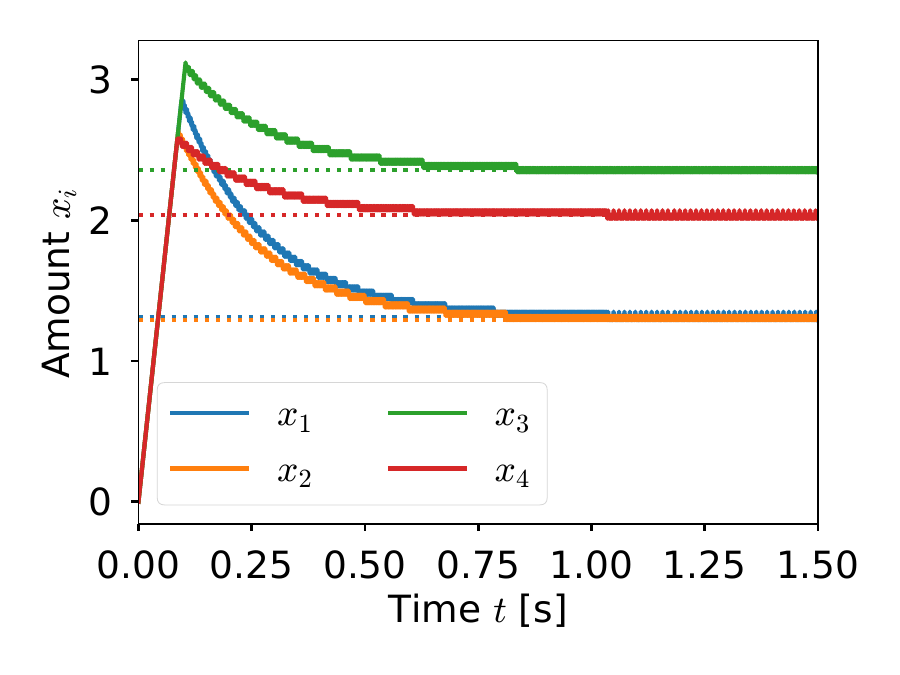}
    \caption{Primal variable $x$ of Sign PDGD with $\eta_x=30.0$}
    \label{F:E3.spdgd-eta30_x}
\end{figure}
\begin{figure}
    \centering
    \includegraphics[width=1\linewidth]{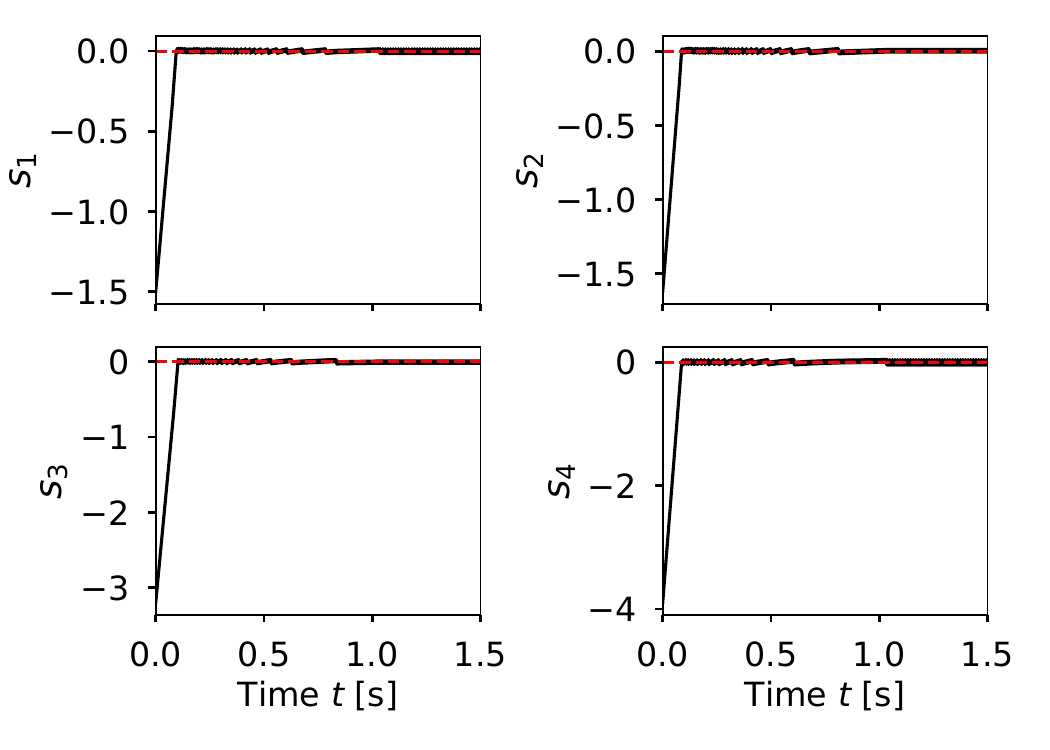}
    \caption{Sliding variable $s$ of Sign PDGD with $\eta_x=30.0$}
    \label{F:E3.spdgd-eta30_s}
\end{figure}

\UpdateA{Fig.~\ref{F:E3.spdgd-eta30_x} shows} the time evolution of $x_i$.
Solid lines show the buyers' trajectories, and the dotted lines show the analytical optimum $x_i^\ast$. All trajectories converge to the analytical optimum.
The allocation of each buyer increases almost linearly up to a certain time. Afterwards, the trajectories bend and settle smoothly to the optimum.

\UpdateA{Fig.~\ref{F:E3.spdgd-eta30_s} shows} the time evolution of the sliding variable $s$. \UpdateA{The black solid line shows the value of $s$ and the dotted red line shows $0$.} Every element first has $s_i\ne0$, then reaches $s_i=0$, after which $s_i=0$ is maintained. This is the behavior predicted for the reaching phase followed by the sliding phase.

\UpdateB{Lemma~\ref{L:A.Reach} gives the upper bound of the reaching phase time as $T^\mathrm{R}=0.136~\text{s}$. The observed time at which the sliding phase started for all elements was $0.123~\text{s}$, which is consistent with the bound\footnote{\UpdateB{Note that the evaluated $\sup_{t\ge0}|\dot v_i(t)|$ is in \eqref{E:A.R.ReachTime.0} from the simulation result, i.e., a posteriori. Hence $T^\mathrm{R}$ cannot be computed before running the dynamics.}}.%
}

\subsection{\UpdateA{Results for $\eta_x=3.0$}}
\begin{figure}
    \centering
    \includegraphics[width=1\linewidth]{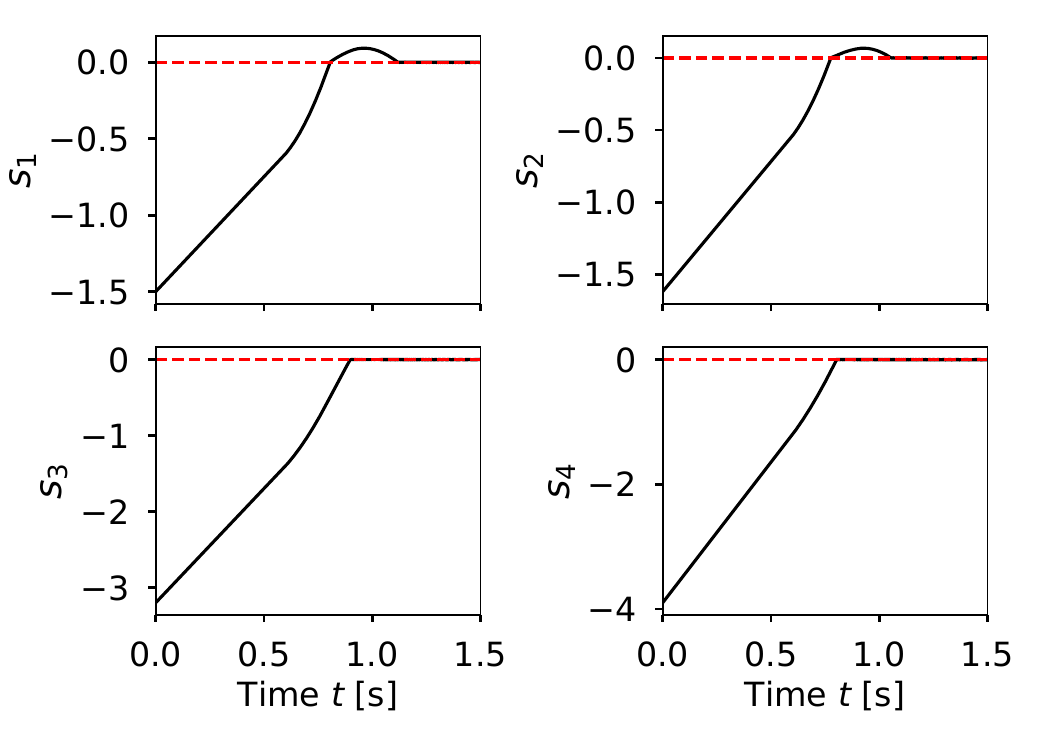}
    \caption{Sliding variable $s$ of Sign PDGD with $\eta_x=3.0$}
    \label{F:E3.spdgd-eta3_s}
\end{figure}

\begin{table*}[t]
    \centering
    \caption{\UpdateA{Check of Assumption~\ref{Assumption:A.PriceVelocityBound}: The maximum price speed and its upper bound}}
    \label{T:E3.PriceVelocityAssumptionCheck}
    \begin{tabular}{c|cccc|cccc}
         \hline
         & \multicolumn{4}{c|}{$\eta_x=30.0$} & \multicolumn{4}{c}{$\eta_x=3.0$} \\
         Element $i$ 
         & 1 & 2 & 3 & 4 
         & 1 & 2 & 3 & 4 \\
         \hline
         $\max_t|\dot v_i(t)|$
             &$3.98$&$3.98$&$3.98$&$3.98$ 
             &$2.49$&$2.49$&$2.49$&$2.49$\\ 
         $\eta_xF_i$
             &$15.00$&$18.00$&$30.00$&$45.00$ 
             &$1.50$&$1.80$&$3.00$&$4.50$\\ 
         \hline
    \end{tabular}
\end{table*}

\UpdateA{%
When $\eta_x=3.0$, the Assumption~\ref{Assumption:A.PriceVelocityBound} does not hold.
Table~\ref{T:E3.PriceVelocityAssumptionCheck} shows the maximum of $|\dot v_i(t)|$ and the bound $\eta_xF_i$ for each element and for both gains $\eta_x$. For $\eta_x=30.0$, the maximum is below the bound for all elements, so the assumption holds.
However, for $\eta_x=3.0$, the maximum exceeds the bound for $i=1$ and $2$, so the assumption is violated for these elements, while it holds for $i=3$ and $4$.

Fig.~\ref{F:E3.spdgd-eta3_s} shows the sliding variable $s$ for $\eta_x=3.0$. 
Elements $i=3$ and $4$, for which Assumption~\ref{Assumption:A.PriceVelocityBound} holds, reach $s_i=0$ and stay there.
}%

\UpdateB{The sliding phase started for all components at $0.898~\text{s}$, which is later than in the $\eta_x = 30.0$ scenario. This indicates that a smaller gain $\eta_x$ slows the reaching phase.}

\UpdateB{Elements $i=1$ and $2$ reach $s_i=0$ but subsequently leave $s_i=0$. After that, they return to $0$ only after a few seconds and converge to the optimum.
The inequality $|\dot v_i(t)| < \eta_x F_i.$ in Assumption~\ref{Assumption:A.PriceVelocityBound} was violated only temporarily. After the violation, the inequality held again, and the sliding phase is maintained. This suggests that Assumption~\ref{Assumption:A.PriceVelocityBound} is sufficient but not necessary, and it may be conservative in this scenario.}

\section{Conclusion}
In this paper, we proposed the Sign PDGD, a primal-dual gradient dynamics for supply-demand coordination in which human buyers provide only coarse feedback.
We modeled the buyer's signal as the sign of the gradient of the buyer utility and marginal utility.
We analyzed the convergence property of the Sign PDGD by interpreting it as SMC.
In the numerical simulation, we assumed a simple resource-sharing problem, and the Sign PDGD converged to the analytical optimum, reproducing the reaching and sliding phase predicted by the theoretical analysis.

Several works remain for future work.
First, a theoretical analysis for the discrete-time system is needed to explain the chattering observed in the numerical simulation.
Second, we assumed that the human intention is given as a sign signal. In practice, a mechanism for extracting this signal from the buyer's behavior is required. For example, \cite{Miyaoka26_ChatMPC} and \cite{Miyaoka26_ChatIaaS} study the use of language models to infer intentions from human language. Combining the intention-reading mechanism with Sign PDGD is a possible direction.
\UpdateB{Third is to derive a weaker sufficient condition than Assumption~\ref{Assumption:A.PriceVelocityBound}.}

\bibliographystyle{IEEETran}
\bibliography{references}

\end{document}